\pdfoutput=1
\documentclass{amsart}
\usepackage{amssymb,amsthm,amsmath,amstext,amsxtra}
\usepackage{bm}       % allows bold italic letters in math mode
\usepackage{mathtools} % allows more extendible arrows
\usepackage[all]{xy}
\usepackage{booktabs}
\usepackage{hyperref}
\hypersetup{colorlinks=true,urlcolor=blue,citecolor=blue,linkcolor=blue}
\usepackage{enumerate}
\usepackage{stmaryrd}
\usepackage{enumitem}
\usepackage{colonequals}
\usepackage{float} % Needed to force the placement of figures
\usepackage[group-minimum-digits=4,group-separator=\text{,}]{siunitx}

\usepackage{tikz}

\newcommand{\tikzcircle}[2][black,fill=black]{\tikz[baseline=-0.5ex]\draw[#1,radius=#2] (0,0) circle ;}%

\usepackage{svg}

\usepackage[color=white, linecolor=red, bordercolor=red]{todonotes}

\usepackage{bold-extra}

\usepackage{wrapfig}
\usepackage[most]{tcolorbox}
\usepackage{graphicx} 
\usepackage{caption}
\usepackage{calc}
\usepackage{subcaption}
\newtheorem{thm}{Theorem}
\newtheorem{prop}[thm]{Proposition}
\newtheorem{lem}[thm]{Lemma}
\theoremstyle{remark}
\newtheorem{rem}[thm]{Remark}

\theoremstyle{definition}

\numberwithin{equation}{section}
\numberwithin{thm}{section}

\usepackage{algpseudocode}
\algtext*{EndWhile}
\algtext*{EndFor}
\algtext*{EndIf}

\algnewcommand{\Forr}[2]{\Statex \hspace{#1}\algorithmicfor \ #2\ \algorithmicdo}
\algnewcommand{\EndForr}{\unskip\ }

\newlist{inputoutputlist}{itemize}{1}
\setlist[inputoutputlist]{label=--,topsep=0pt,leftmargin=20pt}

\numberwithin{equation}{section}

\newcommand{\OO}{\mathcal{O}}

\newcommand{\C}{\mathbb C} % C for complex
\newcommand{\R}{\mathbb R} % R for reals
\newcommand{\Z}{\mathbb Z} % Z for integers
\newcommand{\Q}{\mathbb Q} % Q for rationals
\newcommand{\F}{\mathbb F} % F for some field
\newcommand{\abs}[1]{\lvert #1 \rvert} % absolute values

\usepackage{xcolor}
\definecolor{darkred}{HTML}{CC1F1F}
\definecolor{green}{rgb}{.4,.7,.4}
\definecolor{blue}{rgb}{.2,.6,.75}
\definecolor{pastelb}{HTML}{3333FF}

\definecolor{pastelyellow}{rgb}{0.992157, 0.552941, 0.235294}
\definecolor{pastelorange}{rgb}{0.941176, 0.231373, 0.12549}
\definecolor{pastelred}{rgb}{0.741176, 0., 0.14902}
\definecolor{darkbrown}{rgb}{0.25098, 0., 0.0745098}

\title[Computing zeta functions of algebraic curves]{Computing zeta functions of algebraic curves using Harvey's trace formula}
\author{Madeleine Kyng}
\email{madeleine.kyng@unsw.edu.au}
\address{School of Mathematics and Statistics, University of New South Wales, Sydney NSW 2052, Australia}

\begin{document}

\begin{abstract}
We present a new method for computing the zeta function of an algebraic curve over a finite field. The algorithm relies on a trace formula of Harvey to count points on a plane model of the curve. The zeta function of the curve is then obtained by making corrections at singular points. We report on an implementation and provide some examples in MAGMA which demonstrate an improvement over Tuitman's algorithm.
\end{abstract}

\maketitle

\section{Introduction}

Let $\F_q$ denote the finite field of characteristic $p$ and cardinality $q=p^a$. Let $\widetilde{X}$ be a nonsingular projective curve of genus $g$ over $\F_q$. Recall that the zeta function of $\widetilde{X}$ is defined by
\[Z(\widetilde{X}, T) = \exp\left(\sum\limits_{r=1}^\infty \frac{\lvert \widetilde{X}(\F_{q^r}) \rvert}{r} T^r\right).\]
It follows from the Weil conjectures \cite[Ch.~VIII]{Lorenzini} that $Z(\widetilde{X}, T)$ is a rational function of the form $P(T)/(1-T)(1-qT)$ where $P(T) \in \Z[T]$ is a polynomial of degree $2g$ whose roots all have absolute value $q^{-\frac{1}{2}}$. Furthermore $Z(\widetilde{X},T)$ satisfies the functional equation  \begin{equation}\label{functional eqn}
Z(\widetilde{X}, q^{-1} T^{-1}) = q^{1-g} T^{2-2g} Z(\widetilde{X}, T).
\end{equation}	

\noindent It follows from these facts that to compute $Z(\widetilde{X}, T)$ it suffices to compute $\lvert \widetilde{X}(\F_{q^r})\rvert$ for $r=1,\dots,g$.

The main purpose of this paper is to describe a practical, efficient algorithm for the problem of computing the zeta function of an arbitrary nonsingular curve.

\begin{thm}\label{main theorem}
	There exists an explicit deterministic algorithm with the following properties. The input consists of a prime $p$, a positive integer $a$, a monic irreducible polynomial $\bar{b} \in \F_p[t]$ of degree $a$ defining the finite field $\F_q \cong \F_p[t]/\bar{b}$, and an absolutely irreducible polynomial $\bar{F} \in \F_q[x,y]$ of degree $d \geq 2$. The output is $Z(\widetilde{X}, T)$ where $\widetilde{X}$ is the nonsingular projective curve with function field $\F_q(x)[y]/\langle \bar{F} \rangle$. The algorithm has time complexity \[a^{O(1)} d^{O(1)} p^{\frac{1}{2}+o(1)}.\]
\end{thm}
We do not provide a complete proof of this theorem. We will give an outline of all of the steps involved in the algorithm, from which one could deduce the time complexity of $a^{O(1)} d^{O(1)} p^{\frac{1}{2}+o(1)}$ by referring to the detailed time complexity estimates provided in \cite{Harv} and \cite{bauch2013complexity}.  We have not worked out the exponents in the time complexity stated in Theorem \ref{main theorem}, but intend to present a thorough analysis in a future paper. In Section \ref{implementation and examples} we present concrete examples that demonstrate the utility and generality of our algorithm. 

%There is no known algorithm that computes $Z(\widetilde{X}, T)$ in time polynomial in the size of the input $O(d^2 \log q)$. The most efficient algorithms have a time complexity that is either polynomial in $\log q$ and exponential in $d$, or else exponential in $\log q$ and polynomial in $d$. The existing algorithms of the former kind have only been used in practice for curves of genus $1$ or genus $2$. The existing algorithms of the latter kind have only been used in practice for certain classes of curves. We present in this paper a practical algorithm of the latter kind which is efficient and applicable to any curve. 
Schoof in \cite{Schoof} was the first to demonstrate a deterministic polynomial-time algorithm for computing the zeta function of an arbitrary genus $g=1$ curve $E/\F_q$. The algorithm involves computing the trace of the Frobenius endomorphism modulo a number of small primes $l$, followed by using the Chinese remainder theorem to determine the exact value of $(q+1 - \#E(\F_q))$. Schoof's algorithm, and higher genus variants of it such as \cite{Pila} and \cite{adleman2001counting}, are known today as \textit{$l$-adic algorithms} for computing the zeta function. These $l$-adic algorithms have time complexity polynomial in $\log q$ for fixed $g$, but unfortunately they are badly exponential in $g$.  Extensions of Schoof's algorithm are frequently used for the case of $g=1$, and specialised $l$-adic algorithms for the case of $g=2$ have been useful in practice \cite{gaudry2012genus}, but as yet $l$-adic algorithms have been impractical for the case of an arbitrary curve of genus $g \geq 3$. 

Kedlaya in \cite{Kedlaya} demonstrated an efficient \textit{$p$-adic algorithm} for the problem of computing the zeta function of an  hyperelliptic curve $\widetilde{X}$ over a finite field of odd characteristic. Kedlaya showed how one could apply the machinery of Monsky--Washnitzer cohomology to this problem, computing the zeta function by explicitly computing the action of Frobenius on this $p$-adic cohomology of $\widetilde{X}$. This method proved to be fruitful for extension and generalisation to larger classes of curves, resulting in the development of $p$-adic point-counting algorithms for the case of superelliptic curves \cite{gaudry2001extension}, $C_{ab}$ curves \cite{denef2006counting} and nondegenerate curves \cite{castryck2006computing}.  These Kedlaya-style algorithms have time complexity polynomial in $g$ but exponential in $\log p$, and they are used in practice for curves of genus $g \geq 2$. 

The most general algorithm among the descendants of \cite{Kedlaya} is Tuitman's algorithm \cite{Tuitman}, which can handle \textit{almost all} inputs $\bar{F} \in \F_q[x,y]$. The main drawback of Tuitman's algorithm is that it requires as input a ``good'' characteristic zero lift of the polynomial $\bar{F} \in \F_q[x,y]$. The properties this lift must have are rather technical; they are described in \cite[Ass.~1]{Tuitman}. Such a lift always exists (provided that $p > 2$ and allowing for extension of the base field), but the problem of efficiently computing a good lift for arbitrary $\bar{F} \in \F_q[x,y]$ is difficult. Castryck, Tuitman and Vermeulen have shown in \cite{Tuitman,castryck2018point,castryck2020lifting} how one can compute a suitable lift for inputs $\bar{F}$ that are nondegenerate with respect to their Newton polytope, or that define a curve of geometric genus at most $5$, or that define a curve of arithmetic gonality at most $5$.

There are more general $p$-adic algorithms that can be used to compute the zeta function of an \textit{arbitrary} $n$-dimensional variety over $\F_q$ and that have polynomial time complexity for fixed $p$. Lauder and Wan were the first to demonstrate an algorithm having these properties; they achieved this by using ideas originally presented by Dwork in his proof of the rationality of the zeta function. For fixed $n$, in the case of an $n$-dimensional hypersurface over $\F_q$, Lauder and Wan achieved a time complexity polynomial in $a$, $p$ and the degree $d$ of the defining polynomial \cite[Thm.~37]{lauder_wan}. In \cite{Harv}, Harvey developed an algorithm similar in nature to that of \cite{lauder_wan} but with asymptotically superior time complexity \cite[Thm.~1.2]{Harv}.%; the exponents in the polynomial time bound are smaller for each parameter $a$, $p$, $d$ in \cite[Thm.~1.2]{Harv} than in \cite[Thm.~37]{lauder2006counting}. 

Our algorithm is based on the trace formula for counting points on a hypersurface from \cite[\S3]{Harv}. We present a modified version of this trace formula in Section~\ref{trace formula}. In \cite{Harv}, Harvey suggested that one could compute the zeta function of a nonsingular curve by using the trace formula to count points on a plane model, which gives the right result except possibly at the singularities of that model, and afterwards making corrections for these singularities. We make explicit in Section \ref{section:algorithm} how one does this. 

Let $\bar{F}$ and $\widetilde{X}$ be as in Theorem \ref{main theorem}, and let $X_0$ be the affine curve cut out by $\bar{F}$, with projective closure $X$. By using Harvey's trace formula we can count the points on $X$ in any extension $\F_{q^r}/\F_q$. If $X$ happens to be nonsingular, then we actually have that $X$ is isomorphic to $\widetilde{X}$ over $\F_q$, and by counting points on $X$ in extensions of $\F_q$ of degree up to $g$ we succeed in computing $Z(\widetilde{X}, T)$. 

If $X$ is singular, a naive approach to computing $Z(\widetilde{X}, T)$ using the trace formula would be to compute $Z(X,T)$ and then remove extraneous factors from the numerator, i.e., remove factors whose roots have absolute value $1$ rather than $q^{-\frac{1}{2}}$. One can compute $Z(X,T)$ by bounding the degree of the numerator using Bombieri's bound \cite[Thm. ~1A]{bombieri1978exponential}, and then counting points on $X$ in extensions of $\F_q$ of degree up to this bound. %Bombieri's bound is quadratic in the degree $d$ of $X$, and thus in general is no better than quadratic in the genus $g$ of $\widetilde{X}$. 

Our algorithm does better than the naive approach --- it avoids costly computations of point-counts on $X$ that are performed when using the naive approach in the case that $X$ is singular.  In our algorithm we count only the $\F_{q^r}$-rational points on $X$ for $r=1,\dots,g$, and we determine using the Montes algorithm \cite{Montes} precisely how $\widetilde{X}$ differs from $X$ by factoring ideals related to the singular points of $X$ in maximal orders $\mathcal{O}$ of $\F_q(\widetilde{X})$.  With the point-counts on $X$ for $r=1,\dots,g$ and the extra information about how $\widetilde{X}$ and $X$ differ we determine the values $\abs{\widetilde{X}(\F_{q^r})}$ for $r=1,\dots,g$, and hence compute $Z(\widetilde{X}, T)$. %In the case of $X$ singular we therefore avoid the costly computations of the point-counts on $X$ in extensions of degree larger than $g$ that we would perform using the naive approach; we replace these computations with the relatively cheap computations of the 

 %In the case of $X$ singular this is significantly more efficient than the naive approach because counting points on $X$ in larger extensions is much more costly than factoring the ideals related to the singularities of $X$.   

The main advantage our algorithm has over Tuitman's is that we require no assumptions about the lift $F$ of $\bar{F}$ to characteristic zero. To apply Harvey's trace formula we must lift $\bar{F} \in \F_q[x,y]$ to some $F \in \Z_q[x,y]$, but \textit{any} lift suffices. We have implemented our algorithm in the $q=p$ case in the computer algebra system MAGMA \cite{magma} and made the code publicly available. In Section \ref{implementation and examples} we compare the performance of our implementation against the MAGMA implementation of Tuitman's algorithm. %Our implementation has time complexity with a dependence in $p$ of $p^2$, but by using the deformation reduction techniques in Harvey's paper \cite{Harv} this could be improved to $p^{\frac{1}{2}+o(1)}$.

\section{Harvey's trace formula}\label{trace formula}

In this section, we shall present a generalisation of the trace formula given in Theorem 3.1 of \cite{Harv}. This version allows us to take into account the shape of the polynomial defining a curve (or more generally, hypersurface), and therefore results in a more efficient computation of point-counts than if one were to use a straightforward implementation of the formula from \cite{Harv}. Our version works with the actual Newton polytope of the polynomial, whereas \cite[Thm.~3.1]{Harv} works with a dilation of the standard simplex that contains that Newton polytope. For our point-counting purposes we only need the case of curves, but we shall present the general hypersurface case as it is no harder to state or prove. 

For any domain $R$,  we denote the Laurent polynomial ring $R[x_1, x_1^{-1}, \dots, x_n, x_n^{-1}]$ by $R[x^{\pm}]$, and for $F \in R[x^{\pm}]$ and $u = (u_1, \dots, u_n) \in \Z^n$, we denote by $[F]_u$ the coefficient of $x^u = x_1^{u_1} \cdots x_n^{u_n}$ in $F$. Throughout, we shall use $K$ to denote a convex polytope in $\R^n$ with integral vertices. We will denote by $K_\Z$ the set of integral points in $K$, i.e., $K_\Z := K \cap \Z^n$. For $F \in R[x^{\pm}]$, we denote by $\Delta(F)$ the Newton polytope of $F$, by which we mean the convex hull in $\R^n$ of the finite set $\{u \in \Z^n : [F]_u \neq 0\} \subseteq \Z^n$. We denote by $P_K$ the free $R$-module on the set of monomials with exponents in $K \cap \Z^n$:  \[P_{K} := \bigoplus\limits_{u \in K_\Z} R \, x^u.\] 
For $s \in \Z^+$, let $sK$ denote the $s$-fold dilation of $K$. For two convex polytopes $K_1, K_2$ denote by $K_1 \oplus K_2$ the Minkowski sum $\{v_1 + v_2 : v_1 \in K_1, v_2 \in K_2\}$. One can show that for $F, G \in R[x^{\pm}]$ and $s \in \Z^+$, we have $\Delta(FG) = \Delta(F) \oplus \Delta(G)$ and $\Delta(F^s) = s \Delta(F)$.

For $n \geq 1$, let $\mathbb{A}^n_{\F_q}$ denote affine $n$-space over $\F_q$ with coordinates $x_1, \dots, x_n$. Let $\mathbb{T}_{\F_{q}}^n$ be the affine torus $\{x_1 \cdots x_n \neq 0 \} \subseteq \mathbb{A}^n_{\F_{q}}$. When we refer to a hypersurface $V$ in $\mathbb{T}_{\F_{q}}^n$ cut out by $\bar{F} \in \F_q[x^{\pm}]$, we mean the variety with geometric points \[\{(c_1, \dots, c_n) \in (\overline{\F_q}^*)^n : \bar{F}(c_1, \dots, c_n) = 0\}\] where $\overline{\F_q}$ is an algebraic closure of $\F_q$. We denote by $V(\F_{q^r})$ the set of $\F_{q^r}$-rational points on $V$, i.e., the geometric points \[\{(c_1, \dots, c_n) \in (\F_{q^r}^*)^n : \bar{F}(c_1, \dots, c_n) = 0\}.\]

Let $\Q_q$ be the unique unramified extension of degree $a$ of $\Q_p$, and let $\Z_q$ be the ring of integers of $\Q_q$. Let $\phi : \F_q \to \F_q$ be the absolute Frobenius map $a \mapsto a^p$. We shall also denote by $\phi$ the unique lift of the Frobenius map to a continuous ring automorphism $\Z_q \to \Z_q$. 

We now define maps $\phi, \psi, T_H, A_H$ analogous to the maps $\phi, \psi, T_H, A_H$ in \cite[\S 3.1]{Harv}. Define $\phi, \psi : \Z_q[x^{\pm}] \to \Z_q[x^{\pm}]$ by \[\phi(G) = \sum\limits_{u \in \Z^n} \phi([G]_u) x^{pu}, \hspace{5mm} \psi(G) = \sum\limits_{u \in \Z^n} \phi^{-1}([G]_{pu}) x^u.\]
For $H \in \Z_q[x^{\pm}]$, let $T_H : \Z_q[x^{\pm}] \to \Z_q[x^{\pm}]$ be the multiplication operator $G \mapsto HG$, and let \[A_H = \psi \circ T_{H^{p-1}}.\]
Note that $A_H$ is a $\phi^{-1}$-semilinear map. 
For $m \geq 1$, we also define \[H^{(m)} = (H \cdot \phi(H) \cdots \phi^{m-1}(H))^{p-1}. \]
It is straightforward to check that  \[A_H^m = \psi^m \circ T_{H^{(m)}}\] 
for any $m \geq 1$, and that $A_H$ maps $P_{K}$ into $P_{K}$ for $K \supseteq \Delta(H)$. Note that $A_H^a = \psi^a \circ T_{H^{(a)}}$ is a linear map, it is not just $\phi^{-1}$-semilinear.

The following theorem directly generalises \cite[Thm. ~3.1]{Harv}. 

\begin{thm}\label{trace formula main theorem}
	Let $K \subseteq \R^n$ be a convex polytope with integral vertices, and let $\bar{F} \in \F_q[x^{\pm}]$ be a non-zero Laurent polynomial with $\Delta(\bar{F}) \subseteq K$. Let $V$ be the hypersurface in $\mathbb{T}_{\F_q}^n$ cut out by $\bar{F}$.  Let $F \in \Z_q[x^{\pm}]$ be any lift of $\bar{F}$ with $\Delta(F) \subseteq K$. Let $r, \lambda, \tau$ be positive integers satisfying \[\tau \geq \frac{\lambda}{(p-1)ar}.\] Then \[\abs{V(\F_{q^r})} = (q^r-1)^n \sum\limits_{s = 0}^{\lambda + \tau - 1} \alpha_s \emph{tr}(A_{F^s}^{ar}) \pmod {p^\lambda},\]
	where \[\alpha_s = (-1)^s \sum\limits_{t=0}^{\tau-1} \binom{-\lambda}{t} \binom{\lambda}{s-t} \in \Z,\]
	and where $A_{F^s}^a$ is regarded as a linear operator on the $\Z_q$-module $P_{sK}$. 
\end{thm}

\begin{rem}
	Harvey mistakenly described the map $A_{F^s}$ as a $\Z_q$-linear operator in \cite[Thm. ~3.1]{Harv}. The map $A_{F^s}^a$ is a $\Z_q$-linear operator, but in general $A_{F^s}$ is only $\phi^{-1}$-semilinear. 
\end{rem}

We do not give a proof of this theorem. The proof is identical to the one provided by Harvey for \cite[Thm.~ 3.1]{Harv}, except that we remove all references to degree and replace them with the application of the results $\Delta(G^s) = s \Delta(G)$ and $\Delta(G_1 G_2) = \Delta(G_1) \oplus \Delta(G_2)$. The key idea of Harvey's proof is to construct an indicator function $J$ on $(\Z_{q^r}^*)^n$ that takes the value $1 \in \Z/p^\lambda \Z$ if the input $c \in (\Z_{q^r}^*)^n$ reduces mod $p$ to a zero of $\bar{F}$, and takes the value $0$ otherwise.

As in \cite[Lem.~3.2]{Harv}, we now give a more computationally explicit description of $A_{F^s}^{ar}$. \begin{lem}\label{trace formula matrix}
	Let $F \in \Z_q[x^\pm]$ be as in Theorem \ref{trace formula main theorem}. The matrix of $A_{F^s}^a$ on $P_{sK}$, with respect to the basis $\{x^u : u \in (sK)_\Z\}$, is given by \[\phi^{a-1}(M_s) \cdots \phi(M_s) M_s\] where $M_s$ is the square matrix defined by \[(M_s)_{v,u} = [F^{(p-1)s}]_{pv-u}\]for $u,v \in (sK)_\Z$, and where $\phi$ acts component-wise on matrices. 
\end{lem}

The proof of this lemma is identical to the one provided in \cite{Harv}.  The difference between our Theorem \ref{trace formula main theorem} and Lemma \ref{trace formula matrix} and Harvey's Theorem 3.1 and Lemma 3.2 is that Harvey works with homogeneous $\bar{F}$ in $\F_q[x_0,\dots,x_n]$ rather than Laurent polynomials in $\F_q[x_1, x_1^{-1},\dots,x_n, x_n^{-1}]$. In Harvey's proof of Theorem 3.1, homogenous polynomials act as a convenient book-keeping device --- the extra variable $x_0$ serves as a way to keep track of relations between Newton polytopes. However, given that the trace formula counts points on a hypersurface in an affine torus it is in fact more appropriate to work with Laurent polynomials, and thus our Theorem \ref{trace formula main theorem} is a more natural way to state Harvey's trace formula.

\section{The algorithm}
\label{section:algorithm}

Let $\F_q$, $\bar{F}$ and $\widetilde{X}$ be as in Theorem \ref{main theorem}. Here we shall present an algorithm based on \cite{Harv} that takes $\bar{F}$ as input and outputs $Z(\widetilde{X}, T)$. Our algorithm consists of two independent subroutines. The first subroutine counts the number of $\F_{q^r}$-rational points on the projective plane model $X$ defined by $\bar{F}$ for $r =1, \dots, g$; the second subroutine determines the errors $\lvert \widetilde{X}(\F_{q^r})  \rvert - \lvert X(\F_{q^r}) \rvert$ for each $r$. We refer to the former subroutine as \textsc{CountPlaneModel}, and the latter as \textsc{ComputeCorrections}. We expect the asymptotic time complexity of \textsc{CountPlaneModel} to dominate that of \textsc{ComputeCorrections}.

\subsection{Counting points on a plane model}

Here we describe the \textsc{Count}\textsc{Plane}-\textsc{Model} algorithm. It takes as input a non-zero polynomial $\bar{F} \in \F_q[x,y]$, a positive integer $D$, and a positive integer $\lambda$. The polynomial $\bar{F}$ cuts out an affine plane curve $X_0$, whose projective closure we denote by $X$. The algorithm outputs a list of the point-counts $\lvert X(\F_{q^r}) \rvert \pmod{p^\lambda}$ for $r = 1, \dots, D$.

The idea of the algorithm is to first count the points on $X \cap \mathbb{T}^2$ by applying the trace formula from Theorem \ref{trace formula main theorem} to $\bar{F}$, followed by locating and counting the points on $X \setminus \mathbb{T}^2$. 

Given the input $\bar{F}, \lambda, D$, the first step in \textsc{CountPlaneModel} is to choose an $F \in (\Z_q/p^\lambda \Z_q)[x,y]$ satisfying $F \pmod p = \bar{F}$ and $\Delta(F) = \Delta(\bar{F})$ to which we can apply the trace formula. Here the ring $\Z_q/p^\lambda \Z_q$ is represented as $(\Z/p^\lambda \Z)[t]/b$, where $b \in (\Z/p^\lambda \Z)[t]$ is a monic, degree $a$ lift of the irreducible degree $a$ polynomial $\bar{b} \in \F_p[t]$ defining the extension $\F_q/\F_p$.  In our application of the trace formula, we will ultimately evaluate \begin{equation}\label{trace formula in algorithm}
\lvert (X \cap \mathbb{T}^2)(\F_{q^r}) \rvert \pmod{p^\lambda} = (q^r-1)^2 \sum\limits_{s=0}^{S} \alpha_s \text{tr}(A_{F^s}^{ar}) \pmod{p^\lambda}
\end{equation} 
for all $r=1,\dots,D$. In each of these sums we take $S := \lambda + \tau -1$ with $\tau := \lceil \lambda/(a(p-1)) \rceil$; by Theorem \ref{trace formula main theorem} this choice of $\tau$ guarantees the correct result modulo $p^\lambda$ for every $r \in \Z^+$. We can do these evaluations of \eqref{trace formula in algorithm} efficiently in the manner described below --- this is the same as what is described in \cite[\S3]{Harv}, with a slight modification to account for the fact that we are working with the Newton polygon $\Delta(F)$ of $F$ rather than the polygon $\text{Conv}\{(0,0), (0,d), (d,0)\} \supseteq \Delta(F)$. 

As in \cite[\S3]{Harv}, we iterate $s$ over $1,\dots,S$: for each $s$ we first construct the matrix $M_s$ modulo $p^\lambda$ where $M_s$ is the matrix described in Lemma \ref{trace formula matrix}, then we compute the traces $\text{tr}(A_{F^s}^{ar}) \pmod{p^\lambda}$ for $r=1,\dots,D$. We compute $M_s$ by expanding $F^{(p-1)s}$ and reading off coefficients at monomials $pv-u$ where $u,v \in s \Delta(F) \cap \Z^2$. We compute the traces $\text{tr}(A_{F^s}^{ar}) \pmod{p^\lambda}$ for $r=1,\dots,D$ exactly as in \cite[Lemma~3.4]{Harv}. That is, we compute the matrix of $A_{F^s}^a$ using a modified binary powering algorithm in accordance with Lemma \ref{trace formula matrix}, and we naively compute the required traces by computing $D$ successive powers of that matrix. Having computed all of the required traces $\text{tr}(A_{F^s}^{ar}) \pmod{p^\lambda}$ for $1 \leq r \leq D$ and $1 \leq s \leq S$, we iterate $r$ over $1, \dots, D$ and use \eqref{trace formula in algorithm} to compute $\lvert (X \cap \mathbb{T}^2)(\F_{q^r}) \rvert \pmod{p^\lambda} $ for each $r$. 

Computing $\lvert (X \setminus \mathbb{T}^2)(\F_{q^r}) \rvert \pmod{p^\lambda}$ for $r=1,\dots,D$ can be done by factoring three univariate polynomials over $\F_q$. We count the points on $X$ in $\mathbb{P}^2 \setminus \mathbb{T}^2$ by factoring the polynomials $\bar{f}_0, \bar{f}_1, \bar{f}_2 \in \F_q[t]$ defined by $\bar{f}_0(t) := \bar{F}_{h}(t, 1, 0)$, $\bar{f}_1(t) := \bar{F}_h(0,t,1)$, and $\bar{f}_2(t) := \bar{F}_h(1,0,t)$, where $\bar{F}_h \in \F_q[x_0, x_1, x_2]$ denotes the homogenisation of $\bar{F}$.

 For simplicity of exposition we have described an unoptimised version of the \textsc{CountPlaneModel} algorithm, whose time complexity has a dependence in $p$ of $p^{2+o(1)}$ rather than $p^{\frac{1}{2} + o(1)}$. To achieve the time complexity stated in Theorem \ref{main theorem}, one can use the deformation recurrence technique from \cite[\S4]{Harv} to compute the required matrices $M_s$ rather than using the naive approach of expanding powers of $F$ and reading off the required coefficients. Another optimisation arises from using an individual precision $\lambda_r$ and upper limit $S_r$ for each $r$ rather than the global precision $\lambda$ and upper limit $S = \lambda + \lceil \lambda/(a(p-1)) \rceil - 1$. %If we use individual upper limits $S(r)$ for each $r$ then we need only work with $S = \max\limits_{1 \leq r \leq D}(S(r))$ of the $M_s$ matrices, and we can cut down on the number of matrix multiplications that are performed.

 The most computationally expensive parts of applying the algorithm described above are: (1) the computation of the coefficients of the powers $F^{(p-1)s}$ that are required to construct the matrices $M_s$, and (2) the computation of the matrix powers $A_{F^s}^{ar}$ and their traces. On (2), note that the dimension of the matrix $M_s$ is roughly $s^2 \text{vol}(\Delta(F))$
 as by construction the matrix $M_s$ has dimension equal to the number of integral points in $s \Delta(F)$, and in $\R^2$ the polygon $s \Delta(F)$ has $s^2$ times the volume of $\Delta(F)$.

\subsection{Computing corrections for miscounted points}

Here we describe the \textsc{ComputeCorrections} algorithm. It takes as input an absolutely irreducible polynomial $\bar{F} \in \F_q[x,y]$ of degree $d \geq 2$ and a positive integer $D$, and it outputs a list of the differences $\lvert \widetilde{X}(\F_{q^r}) \rvert - \lvert X(\F_{q^r}) \rvert$ for $r=1,\dots,D$, where $X$ is the projective closure of the affine plane curve cut out by $\bar{F}$. 

The idea of the algorithm is to reduce the problem of computing the differences $\lvert \widetilde{X}(\F_{q^r}) \rvert - \lvert X(\F_{q^r}) \rvert$ to the problem of counting points on $\widetilde{X}$ and $X$ above a certain set of points on $\mathbb{P}^1$. To do this, we rely on the following theorem. \begin{thm}\label{diff thm}
Let $\pi : \widetilde{X} \to X$ be the normalisation of $X$. Let $S$ be the singular subscheme of $X$. Let $Z$ be any subscheme of $X$ containing $S$, and let $\widetilde{Z} = \pi^{-1}(Z)$. 
 Then we have \[\lvert \widetilde{X}(\F_{q^r}) \rvert - \lvert X(\F_{q^r}) \rvert = \lvert \widetilde{Z}(\F_{q^r}) \rvert - \lvert Z(\F_{q^r}) \rvert.\]
\end{thm}
This follows from: (1) the fact that for any closed point $\mathfrak{p}$ on $X$, we have \[\widetilde{\OO_{X, \, \mathfrak{p} }} = \bigcap\limits_{\mathfrak{q} \in \pi^{-1}(\mathfrak{p})} \OO_{\widetilde{X}, \, \mathfrak{q} }\]
where $\OO_{X, \, \mathfrak{p} } \subseteq \F_q(\widetilde{X})$ is the ring of rational functions on $X$ that are regular at $\mathfrak{p}$, $\widetilde{\OO_{X, \, \mathfrak{p} }}$ is the integral closure in $\F_q(\widetilde{X})$ of that ring, and $\OO_{\widetilde{X}, \, \mathfrak{q} } \subseteq \F_q(\widetilde{X})$ is the ring of rational functions on $\widetilde{X}$ that are regular at $\mathfrak{q}$ \cite[Thm.~III.2.6]{stichtenoth2009algebraic}, and (2) the fact that $\OO_{X, \, \mathfrak{p} }$ is integrally closed if $\mathfrak{p}$ is a nonsingular closed point on $X$ \cite[Prop.~VII.2.6]{Lorenzini}. 

We will compute the differences by choosing a subscheme $Z$ of $X$ containing $S$ and counting points on both $Z$ and $\widetilde{Z} = \pi^{-1}(Z)$. Our \textsc{ComputeCorrections} algorithm will consist of three subroutines: \textsc{ComputeY}, \textsc{CountPointsOnZ} and \textsc{CountPointsAboveZ}, which we describe below. 

\subsubsection{Counting points on $Z$} \label{counting points on $Z$}

Let $a_0(x), \dots, a_n(x)$ be the coefficients of $\bar{F}$ when regarded as a polynomial in $y$ over $\F_q[x]$, i.e., when we write $\bar{F} = a_n(x) y^n + \dots + a_0(x)$ with $a_n(x) \neq 0$. Let $X_0$ be the affine plane curve cut out by $\bar{F}$, and let $S_0$ be the singular subscheme of $X_0$. Let $\varphi : X_0 \to \mathbb{A}^1$ be the morphism of curves $(x,y) \mapsto x$. Let $A$ be the closed set of $\mathbb{A}^1$ where $a_n(x)$ vanishes. %Let $\mathbb{A}^2$ be embedded in $\mathbb{P}^2$ by $(x,y) \mapsto [x:y:1]$.

We define subschemes $Y_0, Y$ of $\mathbb{P}^1$ and $Z_0, Z$ of $X$ as follows: 
\begin{equation}\label{subscheme definitions}
	\begin{split}
		Y_0 &:= \varphi(S_0) \cup A, \\
		Z_0 &:= \varphi^{-1}(Y_0), \\
		Y &:= Y_0 \cup \{\infty\}, \\
		Z &:= Z_0 \cup (X \setminus X_0).\\
	\end{split}
\end{equation}
That is, the points of $Y_0$ are the $x$-coordinates of singular points on $X_0$ together with the $x$-coordinates of ``vertical asymptotes'' of $X_0$. The points on $Z_0$ are those points on $X_0$ that share an $x$-coordinate with a singular point, along with those points at which $a_n(x)$ vanishes. The points on $Z \subseteq X$ are all of the points on $Z_0 \subseteq X_0$ together with all of the points on $X \setminus \mathbb{A}^2 \subseteq \mathbb{P}^2$. Note that by construction, $Z$ contains $S$. An illustration	 of the schemes $Z_0$ and $Y_0$  is shown below in Figure \ref{defn_of_scheme}.

\begin{figure}[h]
	\scalebox{1.5}{\includegraphics{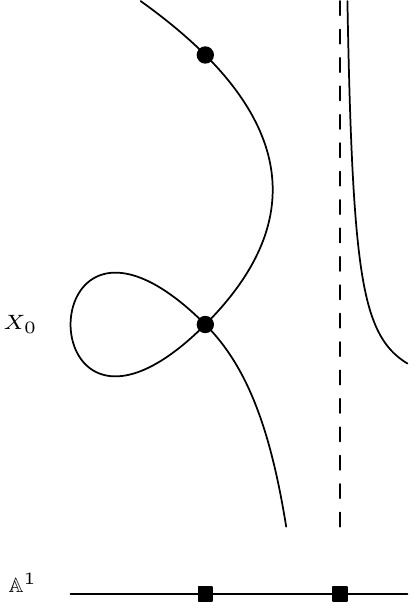}}
	
	%\def\svgwidth{\columnwidth}
	%\scalebox{0.5}{\input{defn_of_scheme_2.pdf_tex}}
	\caption{An example showing $Y_0 \subseteq \mathbb{A}^1$ ($\blacksquare$) and $Z_0 \subseteq X_0$ (\tikzcircle{3.5pt}).}
	\label{defn_of_scheme}
	\centering
\end{figure}

Based on the discussion above we now introduce the subroutine \textsc{ComputeY}. It takes as input the polynomial $\bar{F}$ and outputs a list of irreducible polynomials in $\F_q[x]$ representing the closed set $Y_0 \subseteq \mathbb{A}^1$. This subroutine amounts to finding the distinct irreducible factors of the polynomial $a_n(x) \cdot \gcd(\text{res}_y(\bar{F}, \frac{\partial \bar{F}}{\partial y}), \, \text{res}_y(\bar{F}, \frac{\partial \bar{F}}{\partial x} ) )$. Note that $\text{res}_y(\cdot, \cdot)$ denotes the resultant where the polynomials are regarded as univariate polynomials in $y$ whose coefficients are polynomials in $x$, and hence the set of roots of $\gcd(\text{res}_y(\bar{F}, \frac{\partial \bar{F}}{\partial y}), \, \text{res}_y(\bar{F}, \frac{\partial \bar{F}}{\partial x} ) )\in \F_q[x]$ includes all of the $x$-coordinates of common zeroes of $\bar{F}, \frac{\partial \bar{F}}{\partial x}$ and $\frac{\partial \bar{F}}{\partial y}$, i.e., it includes the $x$-coordinates of singular points on $X_0$.

The algorithm \textsc{CountPointsOnZ} takes as input the polynomial $\bar{F}$, the list of irreducible polynomials representing $Y_0$, and a positive integer $D$, and outputs the list of point-counts $\lvert Z(\F_{q^r}) \rvert$ for $r=1,\dots,D$. Counting points on $Z$ amounts to: (1) counting the points on $Z \cap \mathbb{A}^2 = Z_0$ by computing the distinct factors in the factorisation of $\bar{F}$ modulo each of the irreducibles in $\F_q[x]$ that represent $Y_0$, and (2) counting the points on $Z \setminus \mathbb{A}^2 = X \setminus \mathbb{A}^2$ by computing the distinct factors in the factorisation of the degree $d$ homogeneous part of $\bar{F}$. More precisely, we count points on $Z_0$ by constructing the field $\F_q[x]/\bar{h}(x)$ for each irreducible $\bar{h} \in \F_q[x]$ representing $Y_0$, then factoring the univariate in $y$ polynomial $\bar{F} \pmod{\bar{h}}$ over $\F_q[x]/\bar{h}(x)$; each distinct irreducible factor in such a factorisation corresponds to a $\text{Gal}(\overline{\F_q}/\F_q)$-orbit of points on $X_0$ whose $x$-coordinates are roots of $\bar{h}(x)$.

\subsubsection{Counting points on $\widetilde{Z}$}
Let $\widetilde{\varphi} : \widetilde{X} \to \mathbb{P}^1$ be the morphism of curves defined by the rational function $x \in \F_q(x)[y]/\langle \bar{F} \rangle$. Note that when $\bar{F}$ is absolutely irreducible and of degree $d \geq 2$ this morphism is surjective.  

The following proposition allows us to count points on $\widetilde{Z} = \pi^{-1}(Z)$ by counting the points on $\widetilde{X}$ that lie above the scheme $Y \subseteq \mathbb{P}^1$ defined in ~\eqref{subscheme definitions}.
\begin{prop}\label{theorem x-coordinate}
	
	The subscheme $\widetilde{Z} = \pi^{-1}(Z)$ of $\widetilde{X}$ satisfies \[\widetilde{Z} = \widetilde{\varphi}^{\,-1}(Y).\]

\end{prop}

\noindent The proof of this proposition is a straightforward consequence of the following fact: if $\mathfrak{q} \in \widetilde{X}$ satisfies $\pi(\mathfrak{q}) \in X_0$ then $\widetilde{\varphi}(\mathfrak{q}) = \varphi(\pi(\mathfrak{q}))$, otherwise if $\pi(\mathfrak{q}) \not\in X_0$ then $\widetilde{\varphi}(\mathfrak{q}) \in A \cup \{\infty\}$. The equality $\widetilde{\varphi}(\mathfrak{q}) = \varphi(\pi(\mathfrak{q}))$ for $\pi(\mathfrak{q}) \in X_0$ follows from the fact that the morphisms $\varphi$ and $\widetilde{\varphi}$ are induced by the same rational function $x \in \F_q(x)[y]/\langle \bar{F} \rangle$; this is illustrated below in Figure \ref{figure x-coordinate}. The inclusion $\widetilde{\varphi}(\pi^{-1}(X \setminus X_0)) \subseteq \{\infty\} \cup A$ follows from: (1) if $\frac{1}{x}$ vanishes at $\mathfrak{p} = \pi(\mathfrak{q})$ then $\widetilde{\varphi}(\mathfrak{q}) = \infty$, and (2) if $\frac{1}{x}$ does not vanish at $\mathfrak{p}= \pi(\mathfrak{q})$ but $\frac{1}{y}$ does vanish at $\mathfrak{p}$, then $\frac{1}{y} \in \mathfrak{m}_{\widetilde{X},\, \mathfrak{q}}$, and since $a_n(x) y$ is integral over $\F_q[x]$ it follows that $a_n(x) \in \mathfrak{m}_{\mathbb{P}^1, \, \widetilde{\varphi}(\mathfrak{q})}$. 

\begin{figure}[H]
	\scalebox{2}{\includegraphics{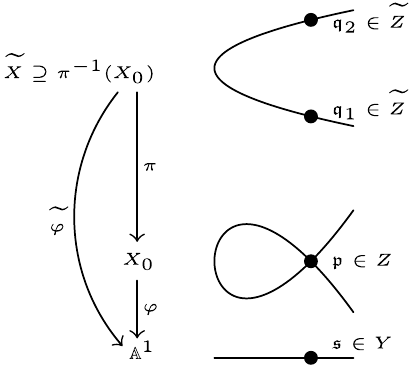}}
	%\def\svgwidth{\columnwidth}
	%\scalebox{0.65}{\input{normalisation_picture_3.pdf_tex}}
	\caption{An example showing $\widetilde{\varphi} = \varphi \circ \pi$.}
	\label{figure x-coordinate}
	\centering
\end{figure}

%Thus in particular the points on $\widetilde{U}$ with $x$-coordinate belonging to $Y_0$ are exactly the points on $\widetilde{U}$ above the points on $U$ with $x$-coordinate belonging to $Y_0$. By our choice of schemes, this means that \[\varphi^{-1}(Y_0) = \pi^{-1}(Z_0) = \pi^{-1}(Z \cap U) = \widetilde{Z} \cap \widetilde{U}.\]

%It is now clear that we can indeed count points on $\widetilde{Z}$ by counting the points on $\widetilde{X}$ lying above $Y$. 
To count the points on $\widetilde{X}$ lying above $Y$, we rely on the Montes algorithm \cite{Montes}, which we shall refer to as \textsc{Montes}. For details on the time complexity of \textsc{Montes} see \cite{bauch2013complexity}. \textsc{Montes} takes as input a monic irreducible separable polynomial $\bar{F} \in \F_q[x][y]$ defining the function field $\F_q(x)[y]/\langle \bar{F} \rangle$ and an irreducible polynomial $\bar{h}(x) \in \F_q[x]$, and outputs the list of closed points on $\widetilde{X}$ that lie above the closed point defined by $\bar{h}(x)$ on $\mathbb{A}^1 = \text{Spec}(\F_q[x])$.

Note that in Theorem \ref{main theorem} we do not assume that the input $\bar{F}$ is monic in $y$, but \textsc{Montes} does require an input which is monic in $y$. Fortunately this is not an issue as it will suffice to apply \textsc{Montes} to the polynomial $\bar{F}' = y^n + a_{n-1}(x) y^{n-1} + \dots + a_0(x) a_n(x)^{n-1}$ instead. This works because $\bar{F}$ and $\bar{F}'$ define the same extension of $\F_q(x)$ via the isomorphism $y \mapsto a_n(x) y$. 

The algorithm \textsc{CountPointsAboveZ} takes as input the polynomial $\bar{F}$, the list of irreducible polynomials representing $Y_0$ and a positive integer $D$, and outputs the list of point-counts $\lvert \widetilde{Z}(\F_{q^r}) \rvert$ for $r=1,\dots,D$. Due to Proposition \ref{theorem x-coordinate}, counting points on $\widetilde{Z}$ amounts to: (1) counting the points on $\widetilde{X}$ lying above $Y_0$ by applying \textsc{Montes} to $\bar{F}'$ for each of the irreducibles $\bar{h}(x)$ in $\F_q[x]$ that represent $Y_0$,  and (2) counting the points on $\widetilde{X}$ that lie above $\infty \in \mathbb{P}^1$  by applying \textsc{Montes} to a polynomial $\bar{F}''$ and $\bar{h}(x) = x$, where $\bar{F}''$ is the polynomial that we get from $\bar{F}'$ via the change of variables $(x,y) \mapsto (\frac{1}{x}, \frac{y}{x^{m}})$ with $m := \deg_x(\bar{F}')$.

\subsubsection{Computing the corrections}

We can now describe the algorithm \textsc{ComputeCorrections}. Given the input $\bar{F} \in \F_q[x,y]$ and $D$, the algorithm first runs the subroutine \textsc{ComputeY} to obtain a representation $L$ of the points in $Y \subseteq \mathbb{P}^1$. It then runs the subroutines \textsc{CountPointsOnZ} and \textsc{CountPointsAboveZ} with arguments $\bar{F}, D$ and $L$ to obtain the point-counts $\lvert Z(\F_{q^r}) \rvert$ and $\lvert \widetilde{Z}(\F_{q^r}) \rvert$ for $r=1,\dots,D$, from which we obtain the differences $\lvert \widetilde{X}(\F_{q^r}) \rvert - \lvert X(\F_{q^r}) \rvert$ for $r=1,\dots,D$ by computing $\lvert \widetilde{X}(\F_{q^r}) \rvert - \lvert X(\F_{q^r}) \rvert = \lvert \widetilde{Z}(\F_{q^r}) \rvert - \lvert Z(\F_{q^r}) \rvert$ in accordance with Theorem \ref{diff thm}.

\subsection{The main algorithm}

We may now outline an efficient algorithm for computing $Z(\widetilde{X}, T)$. Given an absolutely irreducible $\bar{F} \in \F_q[x,y]$ of degree $d \geq 2$, we can compute the genus $g$ of $\F_q(x)[y]/\langle \bar{F} \rangle$ using existing efficient algorithms \cite{hess2002computing}. From the Weil conjectures, we have the Hasse--Weil bound \begin{equation}\label{hasse-weil bound}
	\big\lvert q^{r}+1 - \abs{\widetilde{X}(\F_{q^r})}\big\rvert \leq 2g q^{r/2}, 
\end{equation}
hence we can recover the value $\abs{\widetilde{X}(\F_{q^r})}$ from its reduction modulo any integer strictly larger than $4gq^{r/2}$. Based on this, we take $\lambda = \lceil \log_p(4g q^{g/2} + 1) \rceil$, and run \textsc{CountPlaneModel} with arguments $\bar{F}$, $\lambda$ and $g$. We get a list of point-counts $\abs{X(\F_{q^r})} \pmod {p^\lambda}$ for $r=1,\dots,g$ as output. We then run \textsc{ComputeCorrections} with arguments $\bar{F}$ and $g$ to obtain a list of the differences $\abs{\widetilde{X}(\F_{q^r})} - \abs{X(\F_{q^r})}$ for $r=1,\dots,g$. Adding these lists together modulo $p^\lambda$ gives us the point-counts $\abs{\widetilde{X}(\F_{q^r})} \pmod {p^\lambda}$ for $r=1,\dots,g$. We then recover the exact values $\abs{\widetilde{X}(\F_{q^r})} \in \Z$ using \eqref{hasse-weil bound}. From the power series expansion of \[(1-T)(1-qT) \exp\left(\sum\limits_{r=1}^g \frac{\abs{\widetilde{X}(\F_{q^r})}}{r} T^r \right)\] we obtain the coefficients of the numerator of $Z(\widetilde{X}, T)$ up to the monomial $T^g$. We obtain the remaining coefficients by using the functional equation ~\eqref{functional eqn} for $Z(\widetilde{X}, T)$.

\section{Implementation and examples}\label{implementation and examples}

We have implemented the algorithm described in this paper for prime fields $\F_p$ in the computer algebra system MAGMA \cite{magma}. We have implemented a version with time complexity $p^{2 + o(1)} d^{O(1)}$.
Note that this differs from the time complexity given in Theorem \ref{main theorem} ---  the dependence in $p$ here is $p^2$ instead of $p^{1/2}$. The reason for this difference is that the current implementation involves computing the polynomials $F^{(p-1)s}$ mentioned in Lemma \ref{trace formula matrix} instead of using the deformation-recurrence technique from \cite{Harv} to obtain only the coefficients that we need. In the future we hope to make available an updated implementation which incorporates the deformation-recurrence technique and which works for arbitrary finite fields $\F_q$.

In this section we present and discuss data we have collected about the implementation of our algorithm and the MAGMA implementation of Tuitman's algorithm.  In Section \ref{comparison examples} we provide examples comparing the runtime of our code to that of the MAGMA implementation of Tuitman's algorithm, and in Section \ref{additional examples} we present some data on example computations of our code on inputs that cannot be readily dealt with by existing methods.

Our implementation takes as input an absolutely irreducible $\bar{F} \in \F_p[x,y]$. Tuitman's implementation takes as input a pair $(F,p)$ where $F \in K[x,y]$, $F$ is monic in $y$, $p$ is a prime, and $K$ is a number field in which $p$ is inert. The polynomial $F$ fed into Tuitman's code must define a ``good'' lift of the curve defined by $F \pmod p$ over the field $\Z_K/p \Z_K$, i.e., $F$ must define a lift satisfying \cite[Ass.~1]{Tuitman}. Tuitman's algorithm has time complexity $\widetilde{O}(p d_x^6 d_y^4 a^3 )$,
where $a = [\Z_K / p \Z_K : \Z/p\Z]$, $d_y$ is the degree of $F$ in $x$, and $d_x$ is the degree of $F$ in $y$. 

As mentioned in the introduction, the main advantage of our algorithm over Tuitman's is that it does not require the computation of a ``good'' lift of an input $\bar{F} \in \F_q[x,y]$ to characteristic zero. Our code is capable of handling arbitrary absolutely irreducible inputs $\bar{F} \in \F_p[x,y]$, including those for which one cannot readily compute a lift $F \in \Q[x,y]$ of $\bar{F}$ satisfying Assumption 1 of \cite{Tuitman}. Furthermore, even when it is feasible to compute a good lift of an input $\bar{F}$, the lift one obtains may have properties that cause Tuitman's code to run slowly when compared to the execution of our code on the original $\bar{F}$.

The inputs $\bar{F}$ for which our code is most likely to outperform Tuitman's are those where $\bar{F}$ does not meet Baker's bound on the genus \cite[Thm.~2.4]{beelen_baker}, i.e., inputs $\bar{F}$ for which the geometric genus of the curve defined by $\bar{F}$ is strictly less than the number of integral interior points of $\Delta(\bar{F})$. Note that these are precisely those $\bar{F}$ that define a singular curve in the projective toric surface associated with $\Delta(\bar{F})$. When $\bar{F}$ meets Baker's bound, one can almost always take a naive Newton polygon preserving lift of $\bar{F}$ to obtain a good lift for Tuitman's algorithm \cite[\S2.2]{castryck2018point}. 

When $\bar{F}$ does not meet Baker's bound, determining a lift satisfying \cite[Ass. 1]{Tuitman} is a difficult problem. In some cases this problem can be dealt with by using the methods of \cite{castryck2018point} or \cite{castryck2020lifting} --- these methods apply to inputs $\bar{F}$ that define curves of geometric genus at most $5$ or whose Newton polygon $\Delta(\bar{F})	$  has lattice width at most $5$. The \textit{width} of a convex polytope $K \subseteq \R^n$ along a direction $\mathbf{d} \in \mathbb{R}^n$ is defined to be \[w_{\mathbf{d}}(K) = \max\limits_{\mathbf{x} \in K} \mathbf{d} \cdot \mathbf{x} -  \min\limits_{\mathbf{x} \in K} \mathbf{d} \cdot \mathbf{x}.\] The \textit{lattice width} of a convex polytope $K \subseteq \R^n$ with integral vertices is defined to be \[\text{lw}(K) = \min\limits_{ \mathbf{d \in \Z^n \setminus \{\mathbf{0} \}}} \,\, w_{\mathbf{d}}(K).\] 
For $K \subseteq \R^2$ one can equivalently define the lattice width as the minimal height of a horizontal strip inside which $K$ can be mapped by a unimodular transformation.

%\todomaddie{Consider mentioning that we have included code which checks the output of our implementation, checking against naive point-counting or trace formula. }

The timings obtained in this section were obtained using MAGMA V2.25-7 on a computer with an Nvidia RTX 3090 GPU and an Intel Core i9-12900K CPU with 128GB of memory running Ubuntu 20.04.  
These timings show what is achievable in practice when running a straightforward MAGMA implementation of our algorithm on high-end consumer hardware. The MAGMA supported linear algebra on an RTX 3090 GPU allows for very fast matrix multiplication, making it feasible to apply our algorithm to a large range of curves. 

It should be noted that the MAGMA implementation of Tuitman's algorithm may be able to be significantly improved. Improvements that appeared in MAGMA V2.25-7 were made by MAGMA developer Allan Steel after being sent test inputs from the present author that were found during the writing of this paper. An advantage of our algorithm (particularly the subroutine \textsc{CountPlaneModel}) is that it is comparatively simple and straightforward to implement, and its current implementation has the potential to be sped up substantially through simple low-level optimisations. 

The code for the implementation of our algorithm is available at \url{https://github.com/Maddels/zeta_function_}. Along with the implementation of our algorithm we provide code for testing the correctness of the output that our code produces. The functions provided for testing correctness do not \textit{guarantee} that the output is correct --- failing a test indicates that an error has occurred; passing a test gives partial evidence that the output is consistent with the correct result. We test for correctness by: (1) comparing the point-count that is predicted by the output zeta function for an extension $\F_{p^{i}}$ where $i > g$ with the point-count that is computed by applying the trace formula (\ref{trace formula main theorem}) for the extension $\F_{p^i}$, and (2) comparing the point-counts that are predicted by the output zeta function for small extensions $\F_{p^{i}}$ with point-counts that are computed by naive enumeration. These tests were applied for each of the inputs appearing in Section \ref{additional examples}. Correctness was checked for the inputs in Section \ref{comparison examples} by comparing the output obtained by our implementation to the output obtained by the MAGMA implementation of Tuitman's algorithm.

\subsection{Comparison with Tuitman's algorithm} \label{comparison examples}

We now provide examples where we compare the runtime of our implementation with that of the MAGMA implementation of Tuitman's algorithm. More precisely, we compare the runtime of our algorithm on an input $\bar{F}$ against the combined time of computing a ``good'' lift $F$ of $\bar{F}$ followed by running Tuitman's algorithm on input $F$. For the examples below, we further provide a breakdown of the runtime of both implementations into the runtimes of their key steps. For our implementation, the key steps are: \begin{itemize}
	\item [(A1)] Computing the powers $F^{(p-1)s}$ needed for the construction of the matrices $M_s$.
	
	\item [(A2)] Computing the matrix powers $M_s^r$ and their traces. 
	
	\item [(A3)] All other required computations in the algorithm, including computing the corrections to the point-counts for the plane model. 
\end{itemize}
For the MAGMA implementation of Tuitman's algorithm combined with lifting to obtain a valid input for Tuitman's code, the key steps are: 
\begin{itemize}
	\item [(B1)] Computing a good lift $F$ of $\bar{F}$ to feed into Tuitman's algorithm. 
	\item [(B2)] Computing the objects $\Delta, r(x), s(x,y), W^0, W^\infty, G^0, G^\infty$ described in \cite[\S2]{Tuitman}.
	\item [(B3)] Computing the basis of $H^1(\widetilde{X})$.
	\item [(B4)] Computing the Frobenius lift.
	\item [(B5)] Computing reduction matrices.
	\item [(B6)] Computing the Frobenius matrix, from which the zeta function is computed by computing the characteristic polynomial.
\end{itemize}
Note that steps (B3)--(B6) above correspond to steps I--IV described in \cite{Tuitman}. %Step (B2) above turned out to be very computationally expensive in the examples we collected, particularly when the input $\bar{F}$ did not meet Baker's bound and defined a curve of moderate genus\todomaddie{Check this when we finish our experiments, might need to alter this sentence}. 
To perform (B1) for a given $\bar{F}$, we used MAGMA implementations of the methods described in \cite{castryck2018point, castryck2020lifting}, and when those methods were not applicable we instead tried computing a lift whose singular points are lifts of the singular points on the input plane model. That is, for an input $\bar{F}$ to which the methods of \cite{castryck2018point, castryck2020lifting} do not apply, we attempted to compute an $F \in \Z[x,y]$ that: (1) reduces mod $p$ to $\bar{F}$, and (2) defines a plane curve over $\Q$ whose singularities reduce mod $p$ to the singularities of $X$. This method does not guarantee a lift satisfying \cite[Ass.~1]{Tuitman}, but in the case where the input $\bar{F}$ defines a nodal plane curve we often succeed in finding such a lift.
%Note that empirically, the MAGMA implementation of Tuitman's algorithm has running time with a dependence of $O(H^2)$ in the height $H = H(Q)$ of the input $Q \in \Q[x,y]$. 

Below we present the runtimes for a selection of $10$ inputs $\bar{F}$ that have a variety of features. These 10 inputs along with our code for performing step (B1) can be found on the webpage  \url{https://github.com/Maddels/zeta_function_} in the MAGMA files \texttt{runtime\_comparison\_examples.m} and \texttt{lift\_curve.m}. Inputs 1 and 2 are examples that were made available by Tuitman at the webpage \url{https://github.com/jtuitman/pcc/blob/master/pcc_p/example_p.m},  inputs 3 and 4 were the examples explicitly given in \cite[\S6]{castryck2020lifting} for $d=4$ and $d=5$, input $5$ was generated by specifying its Newton polygon and singular behaviour, and the remaining $5$ inputs were generated using the MAGMA function \texttt{RandomCurveByGenus}. 

Table 1 gives information on the nature of these inputs $\bar{F}$ and the curves they define, as well as which lifting strategy was applied for each input. In Table 1 we denote a naive lift by \texttt{N}, a lift obtained using the methods of \cite{castryck2018point} by \texttt{CT}, a lift obtained using the methods of \cite{castryck2020lifting} by \texttt{CV}, and a singular-point-preserving lift by \texttt{S}. Table 2 gives information on the runtime and memory usage of the implementations of each algorithm on each of these inputs, and Table 3 provides a breakdown of these runtimes into steps A1-A3 and B1-B6. If the runtime for an example exceeded 12 hours then the computation was terminated, with information provided on the runtimes of completed key steps up until the point of termination.

%\begin{table}[H]\centering
%	\caption{Inputs for runtime comparison}
%	\begin{tabular}{|c|c|c|c|c|c|c|c|}
%		\toprule
%		curve & prime  & degree & genus & singular plane model & $\#(\text{int}(\Delta\bar{F}) \cap \Z^2)$  & $\text{lw}(\Delta \bar{F})$&  $\#(\Delta\bar{F} \cap \Z^2)$\\
%		\midrule
%		1 & 7  & 4  & 3 & false & 3 & 4 & 15 \\
%		2 & 11 & 9 & 16 & true & 16 &  5 &  31 \\
%		3 & 17  & 50  & 9 & true & 96 & 5 & 156 \\
%		4 & 23  & 6  & 6 & true & 10 & 6 & 28\\
%		5 & 3 & 8  & 8 & true & 21 & 7 & 41\\
%		6 & 19  & 9  & 10 & true & 28 & 9 & 55 \\
%		7 & 61  & 60  & 5 & true & 145 & 6 & 217 \\
%		8 & 97  & 42  & 5 & true & 100 & 6 & 154\\
%		9 & 11  & 13  & 3 & true & 12 & 3 & 31\\
%		10 & 31  & 10  & 4 & true & 16 & 5 & 36\\
%		\bottomrule
%	\end{tabular}
%	\label{tab:singular points}
%\end{table}

 \begin{table}[h]
	\caption{Inputs for runtime comparison}
	\begin{tabular}{l|cc|ccc|ccc|cc}
	input \# & 1 & 2 & 3 & 4 & 5 & 6 & 7 & 8 & 9 & 10 \\
	\hline
	lifting strategy & \texttt{N} & \texttt{N} & \texttt{CV} & \texttt{CV} & \texttt{CV} & \texttt{S} & \texttt{S} & \texttt{S} & \texttt{CT} & \texttt{CT} \\
	\hline
	meets Bakers bound & yes & yes & no & no & no & no & no & no & no & no \\
	\hline 
	&  &  &  &  &  &  &  &  &  &  \\[-1em]
	$\text{lw}(\Delta(\bar{F}))$ & 3&  5 & 4 & 5 & 5 & 6 & 7 & 9 & 5 & 5 \\
	\hline
	$p$ & $1009$ & $11$ & $7$ & $17$ & $101$ & $23$ & $3$ & $19$ & $43$ & $11$ \\
	\hline
	$g$ & $4$ & $16$ & $10$ & $9$ & $12$ & $6$ & $8$ & $10$ & $4$ & $5$ \\
	\hline 
	&  &  &  &  &  &  &  &  &  &  \\[-1em]
	$\#\big(\text{int}(\Delta\bar{F}) \cap \Z^2 \big) -g$ & $0$ & $0$ & $47$ & $87$ & $54$ &  $4$ & $13$ & $18$ & $2$ & $5$\\
	\hline
	&  &  &  &  &  &  &  &  &  &  \\[-1em]
	$\text{Vol}(\Delta(\bar{F}))$ & $7.5$ & $22.5$ & $80$ & $125$ & $87.5$  & $18$ & $30$ & $40.5$ & $12.5$ & $17.5$ \\
	\end{tabular}
\end{table}

%\todomaddie{Would it be worth including a little (scaled, so each the same size) picture of the Newton polygon of each $\bar{F}$? Or is this overkill?}

%\addtocounter{table}{-1} 

 \begin{table}[h]
	\caption{Runtime and memory comparison}
	\begin{tabular}{c|rr|rr} \label{GPU_compare}
		curve & HK time & T time &   HK memory & T memory \\
		\hline
		$1$ & \texttt{37.98s} & \texttt{\textbf{7.65s}} & \texttt{4141MB} & \texttt{122MB}  \\
		$2$ & \texttt{27.38s} & \texttt{\textbf{14.32s}} & \texttt{1608MB} & \texttt{469MB} \\
		$3$ & \texttt{13.67s} & \texttt{\textbf{1.96s}} & \texttt{922MB} & \texttt{479MB} \\
		$4$ & \texttt{\textbf{40.95s}} & \texttt{9466.71s} & \texttt{2793MB} & \texttt{735MB}  \\
		$5$ & \texttt{\textbf{199.16s}} & \texttt{$\geq$12h} &\texttt{8751MB} & \texttt{$\geq$4415MB}  \\
		$6$ & \texttt{\textbf{0.32s}} & \texttt{437.31s} & \texttt{1835MB} & \texttt{1954MB}  \\
		$7$ & \texttt{\textbf{3.43s}} & \texttt{10607.16s} & \texttt{2032MB} & \texttt{6567MB}\\
		$8$ & \texttt{\textbf{9.55s}} & \texttt{$\geq$12h} & \texttt{2889MB} & \texttt{$\geq$3844MB}\\
		$9$ & \texttt{\textbf{0.08s}} & \texttt{2.01s} & \texttt{2435MB} & \texttt{2435MB}\\
		$10$ & \texttt{\textbf{0.09s}} & \texttt{32.55s} & \texttt{2435MB}& \texttt{2435MB}\\
	\end{tabular}
\end{table}

 %\begin{table}[h]
 %	\caption{Runtime and memory comparison (i7-8700, no GPU)}
%	\begin{tabular}{c|rr|rr} \label{no_GPU_compare}
%		\# & HK time & T time & HK memory & T memory \\
%		\hline
%		$1$ & \texttt{\textbf{0.12s}} & \texttt{0.21s} & \texttt{32MB} & \texttt{32MB} \\
%		$2$ & \texttt{163.00s} & \texttt{\textbf{24.21s}} & \texttt{1608MB} & \texttt{469MB} \\
%		%\hline
%		$3$ & \texttt{113.69s} & \texttt{\textbf{3.32s}} & \texttt{1050MB} & \texttt{479MB} \\
%		$4$ & \texttt{\textbf{369.35s}} & \texttt{14853.62s} & \texttt{2362MB} & \texttt{735MB} \\
%		$5$ & \texttt{\textbf{1086.08s}} & \texttt{$\geq$86400.00s} & \texttt{8710MB} & \texttt{$\geq$4636MB} \\
%		%\hline
%		$6$ & \texttt{\textbf{0.97s}} & \texttt{801.16s} & \texttt{64MB} & \texttt{458MB} \\
%		$7$ & \texttt{\textbf{13.66s}} & \texttt{19617.36s} & \texttt{598MB} & \texttt{5057MB} \\
%		$8$ & \texttt{\textbf{54.13s}} & \texttt{$\geq$86400.00s} & \texttt{1475MB} & \texttt{$\geq$3833MB} \\
		%\hline
%		$9$ & \texttt{\textbf{0.27s}} & \texttt{4.12s} & \texttt{32MB} & \texttt{32MB} \\
%		$10$ & \texttt{\textbf{0.18s}} & \texttt{59.49s} & \texttt{32MB} & \texttt{32MB} \\
%	\end{tabular}
%\end{table}

%\addtocounter{table}{-1} 

%\todomaddieinline{Add \textit{lifting time} to Tuitman runtime, currently it is not included - except for 9 and 10.}

\begin{table}[!htb]
	\caption{HK and T runtime breakdown (in seconds)}
		\hspace*{-0.2cm}
		\begin{tabular}{c|ccc|cccccc}\label{no_GPU_compare_breakdown}
		\# & A1 & A2 & A3 & B1 & B2 & B3 & B4 & B5 & B6\\
		\hline
		1 & \texttt{37.44} & \texttt{0.11} & \texttt{0.43} & \texttt{0.00} & \texttt{0.02} & \texttt{0.04} & \texttt{1.00} & \texttt{1.07} & \texttt{5.52} \\
		2 & \texttt{0.46} & \texttt{26.76} & \texttt{0.16} & \texttt{0.00} & \texttt{0.64} & \texttt{2.08} & \texttt{1.42} & \texttt{0.69} & \texttt{9.49} \\
		3 & \texttt{0.24} & \texttt{13.39} & \texttt{0.04} & \texttt{0.18} & \texttt{0.10} & \texttt{0.45} & \texttt{0.13} & \texttt{0.16} & \texttt{0.94} \\
		4 & \texttt{1.13} & \texttt{39.72} & \texttt{0.10} & \texttt{64.07} & \texttt{73.22} & \texttt{9222.68} & \texttt{1.53} & \texttt{0.62}
		& \texttt{104.58} \\
		5 & \texttt{107.50} & \texttt{90.97} & \texttt{0.69} & \texttt{2236.41} & \texttt{1322.85} & \texttt{-} & \texttt{-} & \texttt{-}
		& \texttt{-} \\
		6 & \texttt{0.15} & \texttt{0.16} & \texttt{0.01} & \texttt{0.10} & \texttt{267.41} & \texttt{125.27} & \texttt{0.69} & \texttt{3.70}
		& \texttt{40.14} \\
		7 & \texttt{0.17} & \texttt{3.22} & \texttt{0.04} & \texttt{0.03} & \texttt{3372.37} & \texttt{7140.78} & \texttt{4.82} & \texttt{19.54}
		& \texttt{69.59} \\
		8 & \texttt{1.13} & \texttt{8.38} & \texttt{0.04} & \texttt{0.58} & \texttt{-} & \texttt{-} & \texttt{-} & \texttt{-}
		& \texttt{-} \\
		9 & \texttt{0.07} & \texttt{0.01} & \texttt{0.00} & \texttt{0.15} & \texttt{0.23} & \texttt{0.37} & \texttt{0.21} & \texttt{0.14}
		& \texttt{0.91} \\
		10 & \texttt{0.02} & \texttt{0.05} & \texttt{0.02} & \texttt{31.82} & \texttt{0.10} & \texttt{0.23} & \texttt{0.05} & \texttt{0.05}
		& \texttt{0.30} \\
		\end{tabular}
\end{table}

%The data from the examples above suggests that our implementation has a reasonable chance to outperform Tuitman's on inputs $\bar{F}$ that do not meet Baker's bound. These are precisely the inputs that require complicated lifting procedures to produce a ``good'' lift, often resulting in Tuitman's algorithm having to work with a different plane model over $\F_q$ than the original input. 

%\subsection{Examples beyond Tuitman's reach} 

\newpage

\subsection{Examples presently beyond Tuitman's algorithm} \label{additional examples}

We now provide example computations of our algorithm on inputs $\bar{F}$ that cannot be readily dealt with by existing methods. We selected our examples based on the lifting strategies described in \S 4.1, in the sense that we selected examples for which these methods cannot be readily applied. One important consideration in the selection of our examples was the lattice width of the Newton polygon $\Delta\bar{F}$. We selected inputs $\bar{F}$ for which it is difficult (or impossible) to find a polynomial $\bar{G} \in \F_q[x,y]$ that both satisfies $\text{lw}(\Delta(\bar{G})) \leq 5$ and defines the same nonsingular curve as $\bar{F}$. 

The \textit{arithmetic gonality} of a curve $\widetilde{X}$ over $\F_q$ is the minimal degree of the extension $\F_q(\widetilde{X})/\F_q(\alpha)$ over all $\alpha \in \F_q(\widetilde{X}) \setminus \F_q$. In the case where $g(\widetilde{X}) \geq 1$, the arithmetic gonality of $\widetilde{X}$ coincides with the minimum lattice width $\text{lw}(\Delta(\bar{G}))$ among all $\bar{G} \in \F_q[x,y]$ that define a function field isomorphic to $\F_q(\widetilde{X})$.

We chose examples $\bar{F}$ with $\text{lw}(\Delta(\bar{F})) \geq 6$ and $g(\widetilde{X}) \geq 9$ as these are likely to have arithmetic gonality exceeding $5$, in which case it would be impossible to find a polynomial $\bar{G} \in \F_q[x,y]$ defining the same function field as $\bar{F}$ to which the methods of \cite{castryck2020lifting} apply. The restriction $g(\widetilde{X}) \geq 9$ is based on the statements about gonality given in \cite[\S2.1]{castryck2018point}. The $\overline{\F_q}$-gonality (i.e., geometric gonality) of a curve over $\overline{\F_q}$ of genus $g \geq 2$ must lie in the range $2, \dots, \lceil g/2 \rceil + 1$, and inside the moduli space $\mathcal{M}_g$ of curves of genus $g$ the dimension of the locus of curves having fixed $\overline{\F_q}$-gonality $ \gamma \in \{2,\dots,\lceil g/2 \rceil + 1\}$ is $\min\{2g + 2 \gamma -5, 3g-3\}$ --- here we are assuming as in \cite{castryck2018point} that this result holds in finite characteristic as well as over $\C$. Thus for $g \geq 9$ the dimension of a locus of curves with an $\overline{\F_q}$-gonality of $\gamma \leq 5$ in $\mathcal{M}_g$ is strictly smaller than that of a locus of curves with an $\overline{\F_q}$-gonality of $\gamma \in \{6,\dots, \lceil g/2 \rceil + 1\}$. 

Furthermore, even for those $\bar{F}$ that define a curve $\widetilde{X}$ of genus $g(\widetilde{X}) \geq 9$ with arithmetic gonality at most $5$, if the lattice width of  $\Delta(\bar{F})$ exceeds $5$ then establishing that the curve has smaller gonality is a difficult problem. One can use the algorithm from \cite{schicho2013computational} to compute geometric gonality, but as explained in \cite{gijswijt2020computing} this becomes impractical beyond genus $7$, and even if the geometric gonality is at most $5$ it may be the case that the arithmetic gonality is larger. 

Another important consideration in the selection of examples $\bar{F}$ was the singular nature of the associated plane curve $X$. We chose examples with complicated singular behaviour, in particular those with non-ordinary singularities. For these inputs it should be more difficult to compute a lift $F$ of $\bar{F}$ that defines a curve of equal genus. 

Ultimately, we selected examples $\bar{F} \in \F_p[x,y]$  with the following properties: \begin{itemize}
	\item[$\bullet$] $\text{lw}(\Delta(\bar{F})) \geq 6$,
	\item[$\bullet$] $g(\widetilde{X}) \geq 9$, 
	\item[$\bullet$] $\bar{F}$ does not meet Baker's bound,
	\item[$\bullet$] the plane curve cut out by $\bar{F}$ has non-ordinary singularities,
	\item[$\bullet$] the naive lift $F \in \Z[x,y]$ of $\bar{F}$, obtained by lifting each coefficient $\bar{c} \in \F_p$ of $\bar{F}$ to an integer $0 \leq c < p$, defines a curve over $\Q$ of genus larger than $g(\widetilde{X})$.
\end{itemize}

%\noindent We furthermore ensured that the examples we chose were not ones that, by chance, could be naively lifted to an $F$ satisfying \cite[Ass. ~1]{Tuitman}. 

In table 5 we present the runtimes of our implementation on a selection of 7 inputs having the above properties. All of the inputs below have a Newton polygon of the form $\text{Conv}\{(0,0), (0,w), (v w, 0)\}$ where $w \geq 6$ and $v \geq 2$. These examples were found by random search among polynomials that: (1) have a Newton polygon that is equal to a specified convex polygon of the form $\text{Conv}\{(0,0), (0,w), (v w, 0)\}$, and (2) define a curve with specified singular behaviour (e.g. having a non-ordinary singularity that takes several blow-ups to resolve). These 7 inputs can be found on the webpage  \url{https://github.com/Maddels/zeta_function_} in the MAGMA file \texttt{examples\_beyond\_tuitman.m}.

 \begin{table}[h]
	\caption{Input curves}
	\begin{tabular}{l|ccccccc}
		input \# & 1 & 2 & 3 & 4 & 5 & 6 & 7\\
		\hline
		&  &  &  &  &  &  &   \\[-1em]
		$\text{lw}(\Delta(\bar{F}))$ & $6$ & $6$ & $6$ & $7$ & $7$ & $8$ & $8$ \\
		\hline
		$p$ & $19$  &  $23$ & $31$ & $19$ & $29$ & $13$ & $17$\\
		\hline
		$g$ & $13$& $20$  & $11$ & $11$ & $12$ & $13$ & $10$ \\
		\hline
		&  &  &  &  &  &  &   \\[-1em]
		$\#(\text{int}(\Delta\bar{F}) \cap \Z^2) -g$ & $72$ & $50$  & $74$ & $67$ & $66$ & $92$  & $95$ \\
		\hline
		&  &  &  &  &  &  &   \\[-1em]
		$\text{Vol}(\Delta(\bar{F}))$ &$108$ & $90$  & $108$ & $98$ & $98$ & $128$ & $128$ \\
	\end{tabular}
\end{table}

 \begin{table}[h] 
	\caption{Runtime and memory usage}
	\begin{tabular}{c|r|rrr|r} \label{GPU_beyond}
		curve &  overall time & A1 & A2 & A3 & memory \\
		\hline
		$1$ & \texttt{253.28s} & \texttt{4.91s} & \texttt{245.82s} & \texttt{2.55s} & \texttt{15097MB} \\
		$2$ & \texttt{1944.40s} & \texttt{29.07s} & \texttt{1908.90s} & \texttt{6.43s} & \texttt{62568MB} \\
		$3$ & \texttt{119.13s} & \texttt{10.06s} & \texttt{108.22s} & \texttt{0.85s} & \texttt{16985MB} \\
		$4$ & \texttt{59.65s} & \texttt{2.76s} & \texttt{56.76s} & \texttt{0.13s} & \texttt{12200MB} \\
		$5$ & \texttt{176.61s} & \texttt{10.80s} & \texttt{164.75s} & \texttt{1.06s} & \texttt{18723MB} \\
		$6$ & \texttt{632.26s} & \texttt{3.40s} & \texttt{625.03s} & \texttt{3.83s} & \texttt{37145MB}\\
		$7$ & \texttt{90.44s} & \texttt{2.51s} & \texttt{87.69s} & \texttt{0.24s} & \texttt{14111MB} \\
	\end{tabular}
\end{table}

\newpage

\section*{Acknowledgments}

Many thanks to David Harvey for his excellent supervision as my doctoral supervisor, for his helpful feedback on earlier versions of this paper, and for suggesting that I pursue the problem of extending his point-counting algorithm. Thanks to Wouter Castryck for helpful discussions on gonality, nondegenerate curves, and the applicability of Tuitman's algorithm. Thanks to Allan Steel and John Cannon of the computational algebra group at the University of Sydney for their assistance with MAGMA and for suggested improvements to the \textsc{CountPlaneModel} algorithm.

%\subsection{Features of our implementation} \label{features of our implementation}

%Here we investigate the time complexity of our implementation. We shall see how the runtime of our implementation varies with several parameters, including the size of the prime $p$, the degree of the input $\bar{F}$, the genus of the curve defined by $\bar{F}$ and the size of the Newton polygon $\Delta \bar{F}$.  

\bibliographystyle{amsalpha}
\bibliography{bib.bib}

\providecommand{\bysame}{\leavevmode\hbox to3em{\hrulefill}\thinspace}
\providecommand{\MR}{\relax\ifhmode\unskip\space\fi MR }
\providecommand{\MRhref}[2]{%
  \href{http://www.ams.org/mathscinet-getitem?mr=#1}{#2}
}
\providecommand{\href}[2]{#2}
\begin{thebibliography}{GSvdW20}

\bibitem[AH01]{adleman2001counting}
Leonard~M. Adleman and Ming-Deh Huang, \emph{Counting points on curves and
  abelian varieties over finite fields}, J. Symbolic Comput. \textbf{32}
  (2001), no.~3, 171--189. \MR{1851164}

\bibitem[BCP97]{magma}
Wieb Bosma, John Cannon, and Catherine Playoust, \emph{The {M}agma algebra
  system. {I}. {T}he user language}, vol.~24, 1997, Computational algebra and
  number theory (London, 1993), pp.~235--265. \MR{1484478}

\bibitem[Bee09]{beelen_baker}
Peter Beelen, \emph{A generalization of {B}aker's theorem}, Finite Fields Appl.
  \textbf{15} (2009), no.~5, 558--568. \MR{2554039}

\bibitem[BNS13]{bauch2013complexity}
Jens-Dietrich Bauch, Enric Nart, and Hayden~D. Stainsby, \emph{Complexity of
  {OM} factorizations of polynomials over local fields}, LMS J. Comput. Math.
  \textbf{16} (2013), 139--171. \MR{3081769}

\bibitem[Bom78]{bombieri1978exponential}
E.~Bombieri, \emph{On exponential sums in finite fields. {II}}, Invent. Math.
  \textbf{47} (1978), no.~1, 29--39. \MR{506272}

\bibitem[CDV06]{castryck2006computing}
W.~Castryck, J.~Denef, and F.~Vercauteren, \emph{Computing zeta functions of
  nondegenerate curves}, IMRP Int. Math. Res. Pap. (2006), Art. ID 72017, 57.
  \MR{2268492}

\bibitem[CT18]{castryck2018point}
Wouter Castryck and Jan Tuitman, \emph{Point counting on curves using a
  gonality preserving lift}, Q. J. Math. \textbf{69} (2018), no.~1, 33--74.
  \MR{3771385}

\bibitem[CV20]{castryck2020lifting}
Wouter Castryck and Floris Vermeulen, \emph{Lifting low-gonal curves for use in
  {T}uitman's algorithm}, A{NTS} {XIV}---{P}roceedings of the {F}ourteenth
  {A}lgorithmic {N}umber {T}heory {S}ymposium, Open Book Ser., vol.~4, Math.
  Sci. Publ., Berkeley, CA, 2020, pp.~109--125. \MR{4235109}

\bibitem[DV06]{denef2006counting}
Jan Denef and Frederik Vercauteren, \emph{Counting points on {$C_{ab}$} curves
  using {M}onsky-{W}ashnitzer cohomology}, Finite Fields Appl. \textbf{12}
  (2006), no.~1, 78--102. \MR{2190188}

\bibitem[GG01]{gaudry2001extension}
Pierrick Gaudry and Nicolas G\"{u}rel, \emph{An extension of {K}edlaya's
  point-counting algorithm to superelliptic curves}, Advances in
  cryptology---{ASIACRYPT} 2001 ({G}old {C}oast), Lecture Notes in Comput.
  Sci., vol. 2248, Springer, Berlin, 2001, pp.~480--494. \MR{1934859}

\bibitem[GMN15]{Montes}
Jordi Gu\`ardia, Jes\'{u}s Montes, and Enric Nart, \emph{Higher {N}ewton
  polygons and integral bases}, J. Number Theory \textbf{147} (2015), 549--589.
  \MR{3276340}

\bibitem[GS12]{gaudry2012genus}
Pierrick Gaudry and \'{E}ric Schost, \emph{Genus 2 point counting over prime
  fields}, J. Symbolic Comput. \textbf{47} (2012), no.~4, 368--400.
  \MR{2890878}

\bibitem[GSvdW20]{gijswijt2020computing}
Dion Gijswijt, Harry Smit, and Marieke van~der Wegen, \emph{Computing graph
  gonality is hard}, Discrete Appl. Math. \textbf{287} (2020), 134--149.
  \MR{4139992}

\bibitem[Har15]{Harv}
David Harvey, \emph{Computing zeta functions of arithmetic schemes}, Proc.
  Lond. Math. Soc. (3) \textbf{111} (2015), no.~6, 1379--1401. \MR{3447797}

\bibitem[Hes02]{hess2002computing}
F.~Hess, \emph{Computing {R}iemann-{R}och spaces in algebraic function fields
  and related topics}, J. Symbolic Comput. \textbf{33} (2002), no.~4, 425--445.
  \MR{1890579}

\bibitem[Ked01]{Kedlaya}
Kiran~S. Kedlaya, \emph{Counting points on hyperelliptic curves using
  {M}onsky-{W}ashnitzer cohomology}, J. Ramanujan Math. Soc. \textbf{16}
  (2001), no.~4, 323--338. \MR{1877805}

\bibitem[Lor96]{Lorenzini}
Dino Lorenzini, \emph{An invitation to arithmetic geometry}, Graduate Studies
  in Mathematics, vol.~9, American Mathematical Society, Providence, RI, 1996.
  \MR{1376367}

\bibitem[LW08]{lauder_wan}
Alan G.~B. Lauder and Daqing Wan, \emph{Counting points on varieties over
  finite fields of small characteristic}, Algorithmic number theory: lattices,
  number fields, curves and cryptography, Math. Sci. Res. Inst. Publ., vol.~44,
  Cambridge Univ. Press, Cambridge, 2008, pp.~579--612. \MR{2467558}

\bibitem[Pil90]{Pila}
J.~Pila, \emph{Frobenius maps of abelian varieties and finding roots of unity
  in finite fields}, Math. Comp. \textbf{55} (1990), no.~192, 745--763.
  \MR{1035941}

\bibitem[Sch85]{Schoof}
Ren\'{e} Schoof, \emph{Elliptic curves over finite fields and the computation
  of square roots mod {$p$}}, Math. Comp. \textbf{44} (1985), no.~170,
  483--494. \MR{777280}

\bibitem[SSW13]{schicho2013computational}
J.~Schicho, F.-O. Schreyer, and M.~Weimann, \emph{Computational aspects of
  gonal maps and radical parametrization of curves}, Appl. Algebra Engrg. Comm.
  Comput. \textbf{24} (2013), no.~5, 313--341. \MR{3118610}

\bibitem[Sti09]{stichtenoth2009algebraic}
Henning Stichtenoth, \emph{Algebraic function fields and codes}, second ed.,
  Graduate Texts in Mathematics, vol. 254, Springer-Verlag, Berlin, 2009.
  \MR{2464941}

\bibitem[Tui17]{Tuitman}
Jan Tuitman, \emph{Counting points on curves using a map to {$\bold{P}^1$},
  {II}}, Finite Fields Appl. \textbf{45} (2017), 301--322. \MR{3631366}

\end{thebibliography}
\end{document}